\documentclass[a4paper]{siamart1116}
\usepackage{mathtools}
\usepackage{epstopdf}
\usepackage{color}
\usepackage[utf8]{inputenc}
\usepackage[T1]{fontenc}
\usepackage{algpseudocode}
\usepackage{algorithm}
\usepackage{amsmath}
\usepackage{amssymb}
\usepackage{float}
\usepackage{breqn}
\usepackage{hyperref}
\graphicspath{{images/}}
\usepackage[square,numbers]{natbib}
\usepackage{subfigure}

\usepackage{hyperref}

\title{Sparse approximate matrix-matrix multiplication for density matrix purification with error control\thanks{\today
\funding{This work was supported by the Swedish national strategic e-science research program (eSSENCE).}}}

\author{Anton G. Artemov%
\thanks{Division of Scientific Computing, Department of Information Technology, Uppsala University, Box 337, SE-751 05 Uppsala, Sweden (\email{anton.artemov@it.uu.se}, \email{emanuel.rubensson@it.uu.se}).}
\and
Emanuel H. Rubensson%
\footnotemark[2]
}

\begin{document}
\maketitle

\begin{abstract}
We propose a method for strict error control in sparse approximate matrix-matrix multiplication. The method combines an error bound and a parameter sweep to select an appropriate threshold value. The scheme for error control and the sparse approximate multiplication are implemented using the Chunks and Tasks parallel programming model. We demonstrate the performance of the method in parallel linear scaling electronic structure calculations using density matrix purification with rigorous error control.  
\end{abstract}

Sparse matrix-matrix multiplication is a key operation in linear
scaling electronic structure calculations based on, for example,
Hartree--Fock or Kohn--Sham density functional theory. This operation
has therefore received a lot of attention from method and software
developers in this field~\cite{DBowler12}. This includes the
development of sparse data structures~\cite{challacombe2000general,
  matrix-ssbrh, m-rrs, densmat-ellpack}, approximation techniques
taking advantage of the special properties of the matrices
involved~\cite{pur-pm, rubensson2011bringing, bock2013optimized}, and
different approaches to parallelization~\cite{borvstnik2014sparse,
  rubensson2016locality, dawson2018massively,
  WeberEtAlMidpointMmul2015, conquest-sparsematrix}.  Sparse matrix-matrix multiplication is
used in the construction of the density matrix defined by
\begin{equation}\label{eq:heavi}
  D = \theta(\mu I - F),
\end{equation}
where $\theta(x)$ is the Heaviside function, $\mu$ is the chemical
potential and $F$ is the Fock or Kohn--Sham matrix.  A number of
different methods for the computation of the density matrix, including
minimization and recursive polynomial expansion methods, use a
sequence of matrix-matrix multiplications. Recursive polynomial
expansion methods are also referred to as density matrix purification.
The great performance of these methods for large systems can be
attributed to the decay properties of the density matrix and any
matrices that occur during the course of its computation. In exact
arithmetics, the matrices involved contain many entries of
small magnitude.
In efficient implementations, this is utilized 
by the removal of small matrix entries
%
%
from the matrix representation, meaning that they are treated
as zeros~\cite{bowler2010calculations, vandevondele2012linear, rudberg2018ergo, dawson2018massively}.
A key issue and a common topic for discussion is how
this is done and what implications it has for performance and
accuracy.

In the recursive polynomial expansion methods it is possible to
strictly control the accuracy of the final result if the norm of the
matrix consisting of removed matrix entries can be controlled. This
procedure is formalized in
Algorithm~\ref{alg:purification_with_truncation} where the removal of
small entries is written as the addition of an error matrix $E_i$ in
each iteration. The starting matrix $X_0$ is given by a linear
transformation of the Fock/Kohn--Sham matrix and $f_i,\, i =
1,2,\dots$ are low order polynomials chosen so that their recursive
application tends to the desired step function in~\eqref{eq:heavi}.
It is shown in refs \citenum{m-accPuri, 2019arXiv190912533K} how to
choose the error tolerances $\delta_i,\, i=1,2,\dots$ so that the error
in the density matrix is controlled, i.e.\ so that
$\|D-\widetilde{D}\| < \varepsilon$, where $\widetilde{D}$ is the
computed approximate density matrix. Several methods to select small
matrix entries for removal without violating the condition $\|E\| <
\delta$ have been proposed, making use of different matrix
norms~\cite{rubensson2011bringing, rubensson2009truncation,
  2019arXiv190912533K}.


\begin{algorithm}[h]
  \caption{Purification process \label{alg:purification_with_truncation}}
  \begin{algorithmic}[1]
    \Procedure{Purification}{$X_0$, $f_i$, $\delta_i$,$n$}
    \State $\widetilde{X}_{0} = X_{0} + E_{0},~ \|E_{0}\|_F < \delta_0$ \label{algline:truncation1}    
    \For{$i = 1,\ldots,n$}
    \State $X_{i} = f_i(\widetilde{X}_{i-1})$
    \State $\widetilde{X}_{i} = X_{i} + E_{i},~ \|E_{i}\|_F < \delta_i$ \label{algline:truncation2}
    \EndFor
    \EndProcedure
  \end{algorithmic}
\end{algorithm}  

A drawback with the approach described above, where the truncation is
performed separately from the multiplication, is that the
multiplication may result in a large increase of the number of nonzero
matrix entries~\cite{hine_onetep, assessment_2011}. Often, many of the
just introduced nonzeros have small magnitude and will anyway be
removed in the subsequent truncation.  This means computational
resources are used to compute and temporarily store those matrix
entries for no purpose.
%
%
Several remedies for this issue have been proposed.  For block-sparse
matrix representations it has been proposed to skip submatrix products
of blocks with small norm~\cite{dmm-c, m-rs}. In the cutoff radius
approach all matrix entries that correspond to distances between
nuclei or basis function centers larger than some cutoff radius are
excluded from the representation~\cite{pur-pm,
  bowler2010calculations}. Since, in this case, the nonzero pattern is
known in advance, the product may be computed directly in truncated
form.

In the present work we are particularly interested in sparse matrix
representations that make use of sparse quaternary trees (quadtree) to
store matrices where any identically zero submatrix quadrant is left
out of the representation~\cite{quadtreeWise1984}.
In the quadtree representation a matrix $X$ is either 1) stored in a
data structure used for small enough matrices, or it is 2) flagged as
identically zero, or it is 3) composed of four quadrants, 
\begin{equation}\label{eq:quadtree}
  X = \begin{pmatrix}
    X_{00} & X_{10} \\ X_{01} & X_{11}
  \end{pmatrix},
\end{equation}
each a
matrix recursively represented as a quadtree.
The data structure used for small matrices may be a regular column- or
row-wise dense matrix representation or some sparse matrix
format.
In the matrix product, zero branches in the quadtree are skipped. In
the SpAMM approach~\cite{challacombe2010fast,bock2013optimized} one
also skips submatrix products whose norm is known to be small. Such
skipping is performed at each level of the quadtree, see
Algorithm~\ref{alg:SPAMM_pseudocode}.

\begin{algorithm}[t]
\caption{SpAMM} \label{alg:SPAMM_pseudocode}
\begin{algorithmic}[1]
\Procedure{SpAMM}{$A,B,\tau$}
\If{$\|A\|_F \|B\|_F < \tau $}
\State{\Return $C = 0$}
\EndIf
\If{lowest level}
\State \Return $C = A B$
\EndIf
\For{$i = 0,1$}
	\For{$j = 0,1$} 
        \State { $T_0 =  \textrm{SpAMM}(A_{i,0}, B_{0,j},\tau)$ } 
        \State { $T_1 =  \textrm{SpAMM}(A_{i,1}, B_{1,j},\tau)$ } 
	\State $C_{i,j} = T_0 + T_1$ 
\EndFor	
\EndFor
\State \Return $C$
\EndProcedure 
\end{algorithmic}
\end{algorithm}

The approaches described above alleviate the issue of fill-in but do
not offer error control.  In this letter, we show how fill-in can be
avoided while strictly controlling the error in the product.
We make use of the SpAMM algorithm but add a preceding step to
carefully select the SpAMM tolerance $\tau$ so that the error in the
product $\|\textrm{SpAMM}(A,B,\tau)-AB\|_F < \delta$, for given $\delta$.
This gives us a method for approximate evaluation $\tilde{f}(X)$ of
low order poynomials $f(X)$ such that $\|\tilde{f}(X)-f(X)\|_F<\delta$
for a given predefined tolerance $\delta$, as required to compute the
density matrix with strict error control.

Our method to select the SpAMM tolerance combines an error bound with
a parameter sweep. We will now show how a bound of the SpAMM product
error can be computed for given input matrices $A$ and $B$ and a given
SpAMM tolerance $\tau$.
Let us consider how the multiplication of 2-by-2 block matrices is
performed with SpAMM. Assume that the blocks are enumerated as
in~\eqref{eq:quadtree}.
%
%
Then the product matrix $C = A B$ is given by
\begin{equation} \label{eq:product_exact}
C = \begin{pmatrix}
A_{00} B_{00} + A_{10} B_{01} & A_{00} B_{10} + A_{10} B_{11}  \\
A_{01} B_{00} + A_{11} B_{01} & A_{01} B_{10} + A_{11} B_{11}
\end{pmatrix}.
\end{equation}
Assume that the SpAMM tolerance $\tau_1$ is such that the whole
procedure is not performed, because $\|A\|_F \| B \|_F < \tau_1.$
Then, clearly, the error matrix $E = C = AB$ and $\|E\|_F \leq \| A
\|_F \|B\|_F < \tau_1$. So the error norm is bounded by the product of
the multiplicands' norms.

Suppose that we multiply the same matrices approximately with some
other tolerance $\tau_2$ and that three of the sub-multiplications
$A_{00} B_{00}, A_{10} B_{01}, A_{11} B_{11}$ are skipped because the
product of norms is too small. The result of this operation is the matrix
\begin{equation} \label{eq:product_approx}
\widetilde{C} = \begin{pmatrix}
0 & A_{00} B_{10} + A_{10} B_{11}  \\
A_{01} B_{00} + A_{11} B_{01} & A_{01} B_{10} 
\end{pmatrix}.
\end{equation} Then, the error matrix $E = C - \widetilde{C}$ is

\begin{align} \label{eq:product_approx_error}
E &= \begin{pmatrix}
A_{00} B_{00} + A_{10} B_{01}  & 0 \\ 0 & A_{11} B_{11}
\end{pmatrix} 
\end{align}  and its Frobenius norm can be bounded from above as

\begin{align}  \label{eq:product_approx_sharper_bound}
\| E \|_F &= \left( \| A_{00} B_{00} + A_{10} B_{01} \|^2_F + \| A_{11} B_{11} \|^2_F \right)^{\frac{1}{2}} \notag \\ 
&\leq \left( \left( \| A_{00} \|_F  \| B_{00}\|_F + \|A_{10} \|_F \| B_{01} \|_F \right)^2 \right. \notag \\ 
&+ \left. \|A_{11} \|^2_F \| B_{11} \|^2_F \right)^{\frac{1}{2}}.
\end{align}

The idea of our algorithm to find the optimal SpAMM tolerance is based
on the observation outlined above: each skipped multiplication brings
an error in the product matrix, and this error can be bounded at any
level of the matrix hierarchy. The summation of the errors from the
underlying multiplications can be done as in
\eqref{eq:product_approx_sharper_bound}.
This gives the error bound for given $A$, $B$, and SpAMM tolerance $\tau$.

In Algorithm \ref{alg:compute_errors} we combine the error bound with
a parameter sweep. This algorithm computes a bound of the SpAMM
product error $\|\textrm{SpAMM}(A,B,\tau)-AB\|_F$ for each of $N$
candidates $(\tau_1,\ldots,\tau_N)$ for the SpAMM tolerance.
Once we know an error bound for each $\tau_i$, it is straightforward to pick the right SpAMM tolerance so that the corresponding error is the largest below the tolerance for the product error.

\begin{algorithm*}[h]
  \caption{\textbf{C}ompute \textbf{S}pAMM \textbf{E}rrors (\textbf{CSE}) \label{alg:compute_errors}}
\begin{algorithmic}[1]
  \Procedure{CSE}{$A,B,(\tau_1,\ldots,\tau_N)$}
  \State $Errors = (0,\ldots,0)$	 \Comment{length $N$}
  \If{$\|A\|_F \|B\|_F = 0$} \Return $Errors$  
  \EndIf
  \If{lowest level}
  \For{$k = 1,\ldots,N$}
	\If{ $\|A\|_F \|B\|_F < \tau_k$ } $Errors[k] = \|A\|_F \|B\|_F$  \EndIf
  \EndFor
  \State \Return $Errors$
  \EndIf
  \For{$i = 0,1$}
	\For{$j = 0,1$}
		\State $E_1 = $ \Call{CSE}{$A_{i0}, B_{0j},(\tau_1,\ldots,\tau_N)$}
		\State $E_2 = $ \Call{CSE}{$A_{i1}, B_{1j},(\tau_1,\ldots,\tau_N)$}
		\State $Errors = Errors + (E_1 + E_2)^{ \circ 2}$ \Comment{Hadamard power, vector sum}
	\EndFor
  \EndFor
  \State $Errors = Errors^{ \circ \frac{1}{2}}$
  \State \Return $Errors$
  \EndProcedure
\end{algorithmic}
\end{algorithm*}

We evaluate our method in the context of density matrix purification with
rigorous error control. In this evaluation, we consider two variants
of Algorithm~\ref{alg:purification_with_truncation}. In the first
variant, we use the new approximate evaluation of the matrix
polynomials with error control, but do not perform any subsequent
truncation on the product, see
Algorithm~\ref{alg:purification_with_approx}. In the second variant,
we include also the subsequent truncation, see
Algorithm~\ref{alg:purification_hybrid}.  Note that in all three
algorithms the error in each iteration, measured by
$\|\widetilde{X}_{i}-f_{i}(\widetilde{X}_{i-1})\|_F$, is controlled by $\delta_i$. Choosing
the error tolerances $\{\delta_i\}$ as described in
ref~\citenum{m-accPuri} allows for strict control of the error in the
final density matrix.

\begin{algorithm}[h]
  \caption{Purification process with approximate multiply \label{alg:purification_with_approx}}
\begin{algorithmic}[1]
  \Procedure{Purification}{$X_0$, $f_i$, $\delta_i$,$n$}
  \State $\widetilde{X}_{0} = X_0$
  \For{$i = 1,\ldots,n$}
  \State $\widetilde{X}_{i} = \tilde{f}_i(\widetilde{X}_{i-1}),~ \|\widetilde{X}_{i}-f_{i}(\widetilde{X}_{i-1})\|_F < \delta_i$
  \EndFor
  \EndProcedure
\end{algorithmic}
\end{algorithm} 


\begin{algorithm}[h]
  \caption{Hybrid purification process \label{alg:purification_hybrid}}
\begin{algorithmic}[1]
  \Procedure{Purification}{$X_0$, $f_i$, $\delta_i$,$n$}
  \State $\widetilde{X}_{0} = X_{0} + E_{0},~ \|E_{0}\|_F < \delta_0$
  \For{$i = 1,\ldots,n$}
  \State $\widehat{X}_{i} = \tilde{f}_i(\widetilde{X}_{i-1}),~ \|\widehat{X}_{i}-f_{i}(\widetilde{X}_{i-1})\|_F < \displaystyle\frac{\delta_i}{2}$ 
  \State $\widetilde{X}_{i} =  \widehat{X}_{i} + E_{i},~ \|E_{i}\|_F < \displaystyle\frac{\delta_i}{2}$
  \EndFor
  \EndProcedure
\end{algorithmic}
\end{algorithm}








We implement the algorithms using the Chunks and Tasks parallel
programming model and library \cite{rubensson2014chunks}.  We use the
Chunks and Tasks matrix library \cite{rubensson2016locality,
  Artemov2018parallelization} and the hierarchical block-sparse leaf
level library of ref~\citenum{2019arXiv190608148A}. The matrix
leaf-level size is 2048, whereas the leaf internal block size is 16.

The computations are performed on the Beskow cluster located at the
PDC center at KTH Royal Institute of Technology in Stockolm,
Sweden. The system is a Cray machine equipped with 2060 nodes each
carrying 2 16-core Intel Xeon E5-2698v3 CPUs combined with 64
gigabytes of RAM. The base operation frequency is 2.3 GHz. The
connection between the nodes is the Cray Aries network with the
Dragonfly topology. The code is compiled with the GCC 8.3.0 compiler,
Cray MPICH 7.7.0 and OpenBLAS 0.2.20. The latter is built from sources
with disabled multi-threading. We let a worker process occupy a whole
computational node. The 32 available cores are split into two groups:
31 perform the tasks in parallel if possible, 1 is dedicated for
communication. 

In our evaluation, we perform two purification iterations with a
converged density matrix using each of the three algorithms, Algorithms
\ref{alg:purification_with_truncation},
\ref{alg:purification_with_approx}, and \ref{alg:purification_hybrid}.
We use a density matrix computed using the Ergo software
\cite{rudberg2018ergo} for a water cluster with 7947 molecules with
the 3-21G basis set, which gives a matrix size 71253.
For a given tolerance, the approximate matrix square is computed. Then
the process is repeated with the same tolerance using the approximate
square from the previous stage as input. In the end, we compute the
exact square of that input matrix to assert that the error does not
exceed the tolerance. We use timings from the second iteration only.


For Algorithm~\ref{alg:compute_errors} we generate a set of possible SpAMM
tolerances by $\tau_1 = 1,\, \tau_i=0.9\tau_{i-1},i=2,\dots,350$,
which gives logarithmically equally spaced values between $1$ and
$10^{-16}$.

%
%
We refer to multiplication of matrices and then truncation as \textit{truncmul}, sparse approximate multiplication as \textit{SpAMM} and their combination as \textit{hybrid}. 

\begin{figure*}[h!]
\center{\includegraphics[width=1\linewidth]{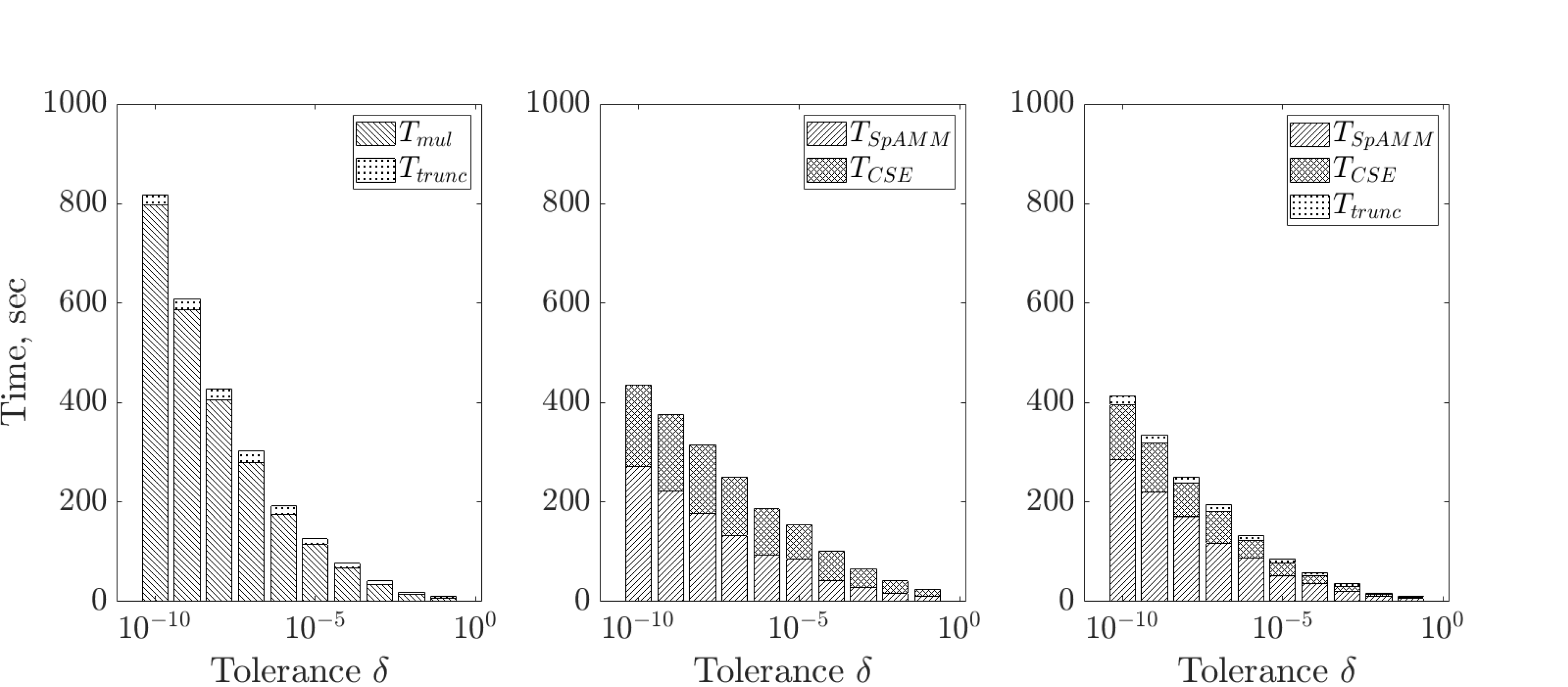}}

\caption{Computational wall time for one iteration in the end of the purification process done with Algorithm \ref{alg:purification_with_truncation} (left), Algorithm \ref{alg:purification_with_approx} (middle) and Algorithm \ref{alg:purification_hybrid} (right), water clusters, matrix size 71253, 16 processes. $T_{mul}$ is the time for multiplication, $T_{trunc}$ is the time for truncation, $T_{SpAMM}$ is the time for the SpAMM Algorithm \ref{alg:SPAMM_pseudocode}, $T_{CSE}$ is the time for the CSE algorithm \ref{alg:compute_errors}.}
\label{fig:waterclusters_timings}
\end{figure*}

We present wall times for the different parts of each of the three
approaches in Figure~\ref{fig:waterclusters_timings}.  The total wall
times of the SpAMM and hybrid approaches with error control proposed
here are less sensitive to the choice of error tolerance and clearly
outperform the truncmul approach for small tolerances.
The hybrid variant outperforms the pure SpAMM variant due to a smaller
time spent on the parameter sweep to select threshold value ($T_{CSE}$).


\begin{figure*}[h!]
\center{\includegraphics[width=1\linewidth]{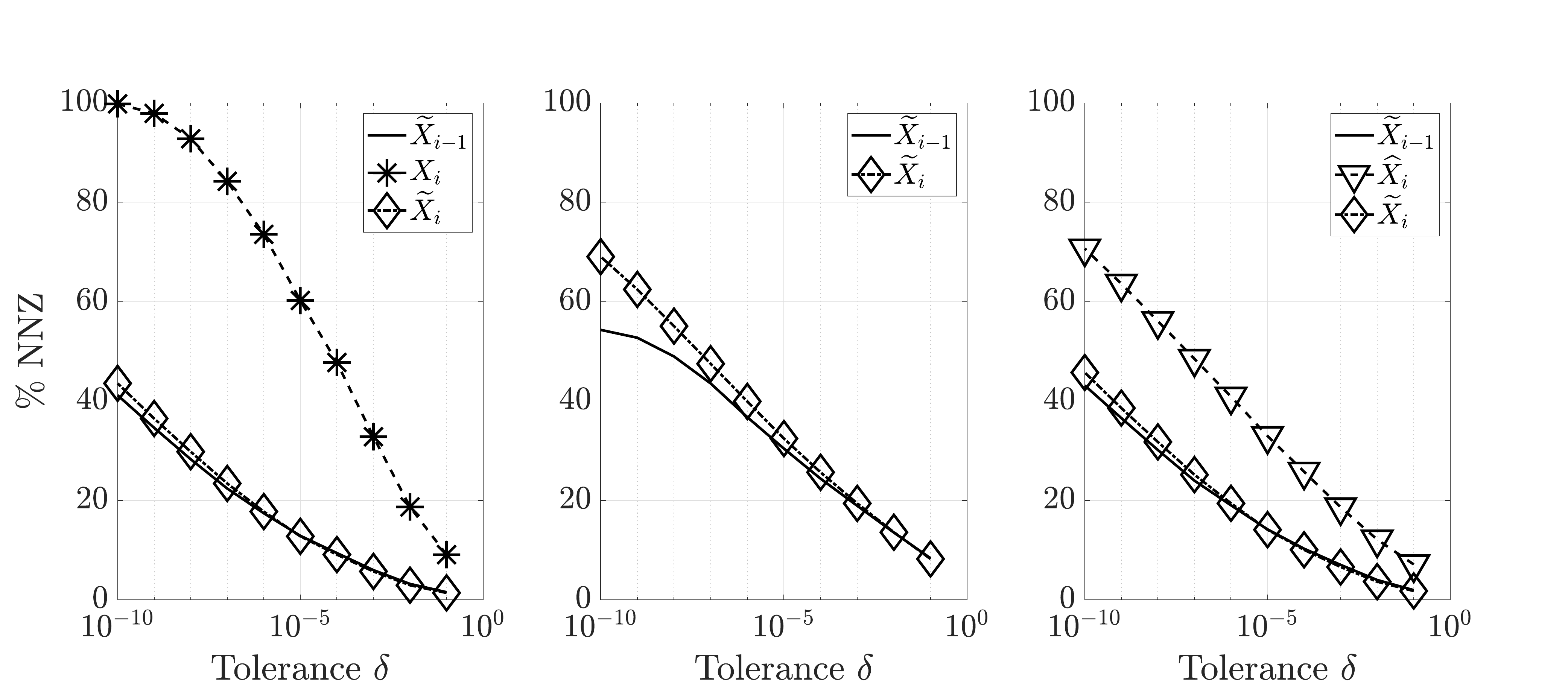}}
\caption{Number of nonzero elements in the end of the purification process done with Algorithm \ref{alg:purification_with_truncation} (left), Algorithm \ref{alg:purification_with_approx} (middle) and Algorithm \ref{alg:purification_hybrid} (right), water clusters, matrix size 71253, 16 processes. Percent of non-zero elements (NNZ) are given for the corresponding matrices.}
\label{fig:waterclusters_nnz}
\end{figure*}

The matrix sparsity for the matrices involved is shown in
Figure~\ref{fig:waterclusters_nnz}.  The left panel clearly shows the
issue discussed earlier with many nonzero entries in $X_i$ computed
for no purpose. Up to around 85\% of the nonzero entries in $X_i$ are
removed in the subsequent truncation. This issue is
mitigated in the SpAMM and hybrid approaches, resulting in smaller
memory usage.


\begin{figure}[h!]
\center{\includegraphics[width=0.75\linewidth]{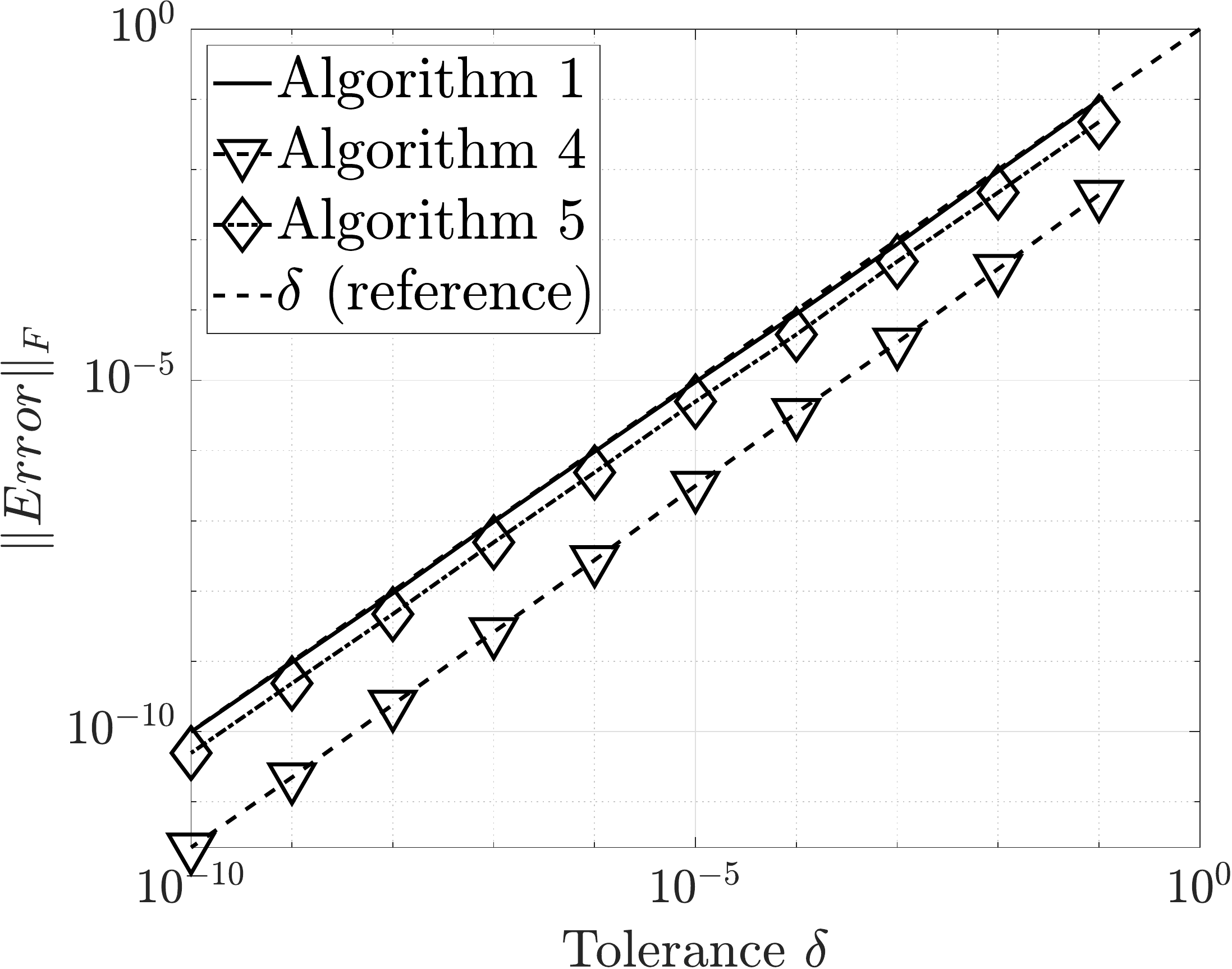}}
\caption{Frobenius norms of error matrices after approximate multiplication of the density matrix, water clusters, matrix size 71253.}
\label{fig:waterclusters_errors}
\end{figure}

We measure the error between the square of the input matrix of the 2nd
iteration of the purification process and its approximate counterpart
computed with each of the approximate multiplication variants to
verify that the error control is working as expected. The results can
be found in Figure \ref{fig:waterclusters_errors}. We can see that all
three variants give an error matrix norm below the desired
tolerance. One can also notice that the truncmul approach provides the sharpest
results in terms of how close the error norm is to the tolerance,
whereas the SpAMM variant provides the least sharp results.

In summary,
we have presented a method to control the Frobenius norm of the error
matrix in sparse approximate matrix-matrix multiplication for matrices
with exponential decay of elements and tested it in the context of the
density matrix purification method. The results show that the standard
routine, see Algorithm \ref{alg:purification_with_truncation}, which
can be described as "multiply-and-truncate" can be improved by
applying the multiplication operation approximately with a properly
chosen threshold. One can build the purification process exclusively
on approximate multiplication, see Algorithm
\ref{alg:purification_with_approx}, or combine it with a subsequent
truncation as done in Algorithm \ref{alg:purification_hybrid}. Our
results indicate that the latter combination is the best option.


Although the new SpAMM and hybrid approaches with error control
clearly outperform the truncmul approach, there is also room
for improvements.  The routines utilizing the SpAMM algorithm require
an extra step, which selects the best truncation threshold value from a
given vector of possible values, and the overhead of this operation is
comparable with the cost of the approximate operation itself for the
variant built exclusively on SpAMM. The hybrid variant has a lower
overhead of the selection routine, which is due to a smaller number of
nonzero elements. The cost of the selection routine
depends on the structure of the matrix, and the more zero blocks it
has, the faster the routine works.
%
%
Another way to reduce the cost is to manipulate the vector of possible
threshold values, for example by altering its length and starting
value.


While representing a significant improvement, the hybrid approach
still involves the computation of a significant number of matrix
entries that are thrown away in the subsequent truncation. This can,
at least partially, be explained by an overestimation of the error by
the CSE algorithm which in Figure~\ref{fig:waterclusters_errors} is
manifested by an error with magnitude more than an order lower than
the chosen tolerance.  Besides improving the error bound in the
Frobenius norm, both with respect to sharpness and speed, as discussed
above, future work could also address error control for SpAMM in other
norms. 

We note that asymptotic error analyses with respect to both the SpAMM tolerance and system size have been carried out in earlier work~\cite{2019arXiv190608148A, 2015arXiv150805856C}. Here, we have proposed a scheme to select the SpAMM tolerance so that the error in a unitary invariant norm is below a predefined tolerance as needed in density matrix purification with rigorous error control. 




\section*{Acknowledgment}
This work was supported by the Swedish national strategic e-science research program (eSSENCE). The computations were performed on resources provided by the Swedish National Infrastructure for Computing (SNIC) at the PDC Center for High Performance Computing, KTH Royal Institute of Technology.

\bibliography{references}

\end{document}